\documentclass[a4paper]{article}
\input epsf.tex
\usepackage{latexsym}
\usepackage{amsmath,amsthm,amsfonts,amssymb,amscd}
\usepackage{enumerate}
\usepackage{amsfonts,amssymb,amsthm,graphicx}
\usepackage[totalwidth=17cm,totalheight=24cm]{geometry}

\catcode`\@=11\relax
\newwrite\@unused
\def\typeout#1{{\let\protect\string\immediate\write\@unused{#1}}}
\def\@nnil{\@nil}
\def\@empty{}
\def\@psdonoop#1\@@#2#3{}
\def\@psdo#1:=#2\do#3{\edef\@psdotmp{#2}\ifx\@psdotmp\@empty \else
    \expandafter\@psdoloop#2,\@nil,\@nil\@@#1{#3}\fi}
\def\@psdoloop#1,#2,#3\@@#4#5{\def#4{#1}\ifx #4\@nnil \else
       #5\def#4{#2}\ifx #4\@nnil \else#5\@ipsdoloop #3\@@#4{#5}\fi\fi}
\def\@ipsdoloop#1,#2\@@#3#4{\def#3{#1}\ifx #3\@nnil 
       \let\@nextwhile=\@psdonoop \else
      #4\relax\let\@nextwhile=\@ipsdoloop\fi\@nextwhile#2\@@#3{#4}}
\def\@tpsdo#1:=#2\do#3{\xdef\@psdotmp{#2}\ifx\@psdotmp\@empty \else
    \@tpsdoloop#2\@nil\@nil\@@#1{#3}\fi}
\def\@tpsdoloop#1#2\@@#3#4{\def#3{#1}\ifx #3\@nnil 
       \let\@nextwhile=\@psdonoop \else
      #4\relax\let\@nextwhile=\@tpsdoloop\fi\@nextwhile#2\@@#3{#4}}
\def\psdraft{
        \def\@psdraft{0}
}
\def\psfull{
        \def\@psdraft{100}
}
\psfull
\newif\if@prologfile
\newif\if@postlogfile
\newif\if@noisy
\def\pssilent{
        \@noisyfalse
}
\def\psnoisy{
        \@noisytrue
}
\psnoisy
\newif\if@bbllx
\newif\if@bblly
\newif\if@bburx
\newif\if@bbury
\newif\if@height
\newif\if@width
\newif\if@rheight
\newif\if@rwidth
\newif\if@clip
\newif\if@verbose
\def\@p@@sclip#1{\@cliptrue}
\def\@p@@sfile#1{
                   \def\@p@sfile{#1}
}
\def\@p@@sfigure#1{\def\@p@sfile{#1}}
\def\@p@@sbbllx#1{
                \@bbllxtrue
                \dimen100=#1
                \edef\@p@sbbllx{\number\dimen100}
}
\def\@p@@sbblly#1{
                \@bbllytrue
                \dimen100=#1
                \edef\@p@sbblly{\number\dimen100}
}
\def\@p@@sbburx#1{
                \@bburxtrue
                \dimen100=#1
                \edef\@p@sbburx{\number\dimen100}
}
\def\@p@@sbbury#1{
                \@bburytrue
                \dimen100=#1
                \edef\@p@sbbury{\number\dimen100}
}
\def\@p@@sheight#1{
                \@heighttrue
                \dimen100=#1
                \edef\@p@sheight{\number\dimen100}
}
\def\@p@@swidth#1{
                \@widthtrue
                \dimen100=#1
                \edef\@p@swidth{\number\dimen100}
}
\def\@p@@srheight#1{
                \@rheighttrue
                \dimen100=#1
                \edef\@p@srheight{\number\dimen100}
}
\def\@p@@srwidth#1{
                \@rwidthtrue
                \dimen100=#1
                \edef\@p@srwidth{\number\dimen100}
}
\def\@p@@ssilent#1{ 
                \@verbosefalse
}
\def\@p@@sprolog#1{\@prologfiletrue\def\@prologfileval{#1}}
\def\@p@@spostlog#1{\@postlogfiletrue\def\@postlogfileval{#1}}
\def\@cs@name#1{\csname #1\endcsname}
\def\@setparms#1=#2,{\@cs@name{@p@@s#1}{#2}}
\def\ps@init@parms{
                \@bbllxfalse \@bbllyfalse
                \@bburxfalse \@bburyfalse
                \@heightfalse \@widthfalse
                \@rheightfalse \@rwidthfalse
                \def\@p@sbbllx{}\def\@p@sbblly{}
                \def\@p@sbburx{}\def\@p@sbbury{}
                \def\@p@sheight{}\def\@p@swidth{}
                \def\@p@srheight{}\def\@p@srwidth{}
                \def\@p@sfile{}
                \def\@p@scost{10}
                \def\@sc{}
                \@prologfilefalse
                \@postlogfilefalse
                \@clipfalse
                \if@noisy
                        \@verbosetrue
                \else
                        \@verbosefalse
                \fi
}
\def\parse@ps@parms#1{
                \@psdo\@psfiga:=#1\do
                   {\expandafter\@setparms\@psfiga,}}
\newif\ifno@bb
\newif\ifnot@eof
\newread\ps@stream
\def\bb@missing{
        \if@verbose{
                \typeout{psfig: searching \@p@sfile \space  for bounding box}
        }\fi
        \openin\ps@stream=\@p@sfile
        \no@bbtrue
        \not@eoftrue
        \catcode`\%=12
        \loop
                \read\ps@stream to \line@in
                \global\toks200=\expandafter{\line@in}
                \ifeof\ps@stream \not@eoffalse \fi
                \@bbtest{\toks200}
                \if@bbmatch\not@eoffalse\expandafter\bb@cull\the\toks200\fi
        \ifnot@eof \repeat
        \catcode`\%=14
}       
\catcode`\%=12
\newif\if@bbmatch
\def\@bbtest#1{\expandafter\@a@\the#1
\long\def\@a@#1
\long\def\bb@cull#1 #2 #3 #4 #5 {
        \dimen100=#2 bp\edef\@p@sbbllx{\number\dimen100}
        \dimen100=#3 bp\edef\@p@sbblly{\number\dimen100}
        \dimen100=#4 bp\edef\@p@sbburx{\number\dimen100}
        \dimen100=#5 bp\edef\@p@sbbury{\number\dimen100}
        \no@bbfalse
}
\catcode`\%=14
\def\compute@bb{
                \no@bbfalse
                \if@bbllx \else \no@bbtrue \fi
                \if@bblly \else \no@bbtrue \fi
                \if@bburx \else \no@bbtrue \fi
                \if@bbury \else \no@bbtrue \fi
                \ifno@bb \bb@missing \fi
                \ifno@bb \typeout{FATAL ERROR: no bb supplied or found}
                        \no-bb-error
                \fi
                \count203=\@p@sbburx
                \count204=\@p@sbbury
                \advance\count203 by -\@p@sbbllx
                \advance\count204 by -\@p@sbblly
                \edef\@bbw{\number\count203}
                \edef\@bbh{\number\count204}
}
\def\in@hundreds#1#2#3{\count240=#2 \count241=#3
                     \count100=\count240        
                     \divide\count100 by \count241
                     \count101=\count100
                     \multiply\count101 by \count241
                     \advance\count240 by -\count101
                     \multiply\count240 by 10
                     \count101=\count240        
                     \divide\count101 by \count241
                     \count102=\count101
                     \multiply\count102 by \count241
                     \advance\count240 by -\count102
                     \multiply\count240 by 10
                     \count102=\count240        
                     \divide\count102 by \count241
                     \count200=#1\count205=0
                     \count201=\count200
                        \multiply\count201 by \count100
                        \advance\count205 by \count201
                     \count201=\count200
                        \divide\count201 by 10
                        \multiply\count201 by \count101
                        \advance\count205 by \count201
                     \count201=\count200
                        \divide\count201 by 100
                        \multiply\count201 by \count102
                        \advance\count205 by \count201
                     \edef\@result{\number\count205}
}
\def\compute@wfromh{
                \in@hundreds{\@p@sheight}{\@bbw}{\@bbh}
                \edef\@p@swidth{\@result}
}
\def\compute@hfromw{
                \in@hundreds{\@p@swidth}{\@bbh}{\@bbw}
                \edef\@p@sheight{\@result}
}
\def\compute@handw{
                \if@height 
                        \if@width
                        \else
                                \compute@wfromh
                        \fi
                \else 
                        \if@width
                                \compute@hfromw
                        \else
                                \edef\@p@sheight{\@bbh}
                                \edef\@p@swidth{\@bbw}
                        \fi
                \fi
}
\def\compute@resv{
                \if@rheight \else \edef\@p@srheight{\@p@sheight} \fi
                \if@rwidth \else \edef\@p@srwidth{\@p@swidth} \fi
}
\def\compute@sizes{
        \compute@bb
        \compute@handw
        \compute@resv
}
\def\psfig#1{\vbox {
        \ps@init@parms
        \parse@ps@parms{#1}
        \compute@sizes
        \ifnum\@p@scost<\@psdraft{
                \if@verbose{
                        \typeout{psfig: including \@p@sfile \space }
                }\fi
                \special{ps::[begin]    \@p@swidth \space \@p@sheight \space
                                \@p@sbbllx \space \@p@sbblly \space
                                \@p@sbburx \space \@p@sbbury \space
                                startTexFig \space }
                \if@clip{
                        \if@verbose{
                                \typeout{(clip)}
                        }\fi
                        \special{ps:: doclip \space }
                }\fi
                \if@prologfile
                    \special{ps: plotfile \@prologfileval \space } \fi
                \special{ps: plotfile \@p@sfile \space }
                \if@postlogfile
                    \special{ps: plotfile \@postlogfileval \space } \fi
                \special{ps::[end] endTexFig \space }
                \vbox to \@p@srheight true sp{
                        \hbox to \@p@srwidth true sp{
                                \hss
                        }
                \vss
                }
        }\else{
                \vbox to \@p@srheight true sp{
                \vss
                        \hbox to \@p@srwidth true sp{
                                \hss
                                \if@verbose{
                                        \@p@sfile
                                }\fi
                                \hss
                        }
                \vss
                }
        }\fi
}}
\catcode`\@=12\relax

\def\E{{\mathbb E}}
\def\pn{{\mathbf p = [\mb{p}_1; \dots; \mb{p}_n]}}
\def\ppn{{\mathbf p' = (\p'_ 1, \dots, \p'_ n)}}
\def\qn{{\mathbf q = [\mb{q}_1; \dots; \mb{q}_n]}}
\def\p{{\mathbf p}}
\def\q{{\mathbf q}}
\def\on{{\mathbf{\omega} = (\dots, \omega_{ij}, \dots)}}
\def\o{{\mathbf{\omega}}}

\newtheorem{prop}{Proposition}[section]
\newtheorem{defi}[prop]{Definition}
\newtheorem{lemma}[prop]{Lemma}
\newtheorem{thm}[prop]{Theorem}
\newtheorem{cor}[prop]{Corollary}

\newtheorem{remark}[prop]{Remark}

\newtheorem{qu}[prop]{Question}

\newcommand{\II}{I\hspace{-0.1cm}I}
\newcommand{\III}{I\hspace{-0.1cm}I\hspace{-0.1cm}I}
\newcommand{\dr}{\partial}
\newcommand{\hess}{\mbox{Hess}}
\newcommand{\tr}{\mbox{tr}}
\newcommand{\be}{\begin{eqnarray}}
\newcommand{\ee}{\end{eqnarray}}
\newcommand{\C}{{\mathbb C}}
\newcommand{\N}{{\mathbb N}}
\newcommand{\R}{{\mathbb R}}
\newcommand{\Z}{{\mathbb Z}}
\newcommand{\mb}{\mathbf}

\newcommand{\eb}{\overline{e}}
\newcommand{\lb}{\overline{l}}
\newcommand{\alphab}{\overline{\alpha}}
\newcommand{\rb}{\overline{r}}

\begin{document}

\title{On the infinitesimal rigidity of weakly convex polyhedra}
\author{
Robert Connelly
\thanks{(visiting Cambridge University until August 2006) Research supported in part by NSF Grant No. DMS-0209595 }\\
Department of Mathematics\\
Malott Hall\\
Cornell University\\
Ithaca, NY 14853\\
USA\\
\texttt{connelly@math.cornell.edu}
\and 
Jean-Marc Schlenker \\
Laboratoire Emile Picard, UMR CNRS 5580\\
Institut de Math{\'e}matiques\\ 
Universit{\'e} Paul Sabatier\\
31062 Toulouse Cedex 9\\
France\\
\texttt{http://www.picard.ups-tlse.fr/\~{ }schlenker}
}

\maketitle

\begin{abstract}
The main motivation here is a question: whether any polyhedron which can be
subdivided into convex pieces without adding a vertex, and which has the same
vertices as a convex polyhedron, is infinitesimally rigid. We prove that it
is indeed the case for two classes of polyhedra: those obtained from a convex
polyhedron by ``denting'' at most two edges at a common vertex, and suspensions with a natural subdivision.
\end{abstract}

\section{A question on the rigidity of polyhedra}

\paragraph{A question.}

The rigidity of Euclidean polyhedra has a long and interesting
history. Legendre \cite{legendre} and Cauchy \cite{Cauchy} proved
that convex polyhedra are rigid: if there is a continuous map between the
surfaces of  two convex polyhedra that is a congruence when restricted to each
face, then the  
map is a congruence between the polyhedra (see \cite{sabitov-legendre}).
However the rigidity of non-convex polyhedra remained an open 
question until the first example of flexible (non-convex) polyhedra 
were discovered \cite{connelly}. 

We say that a polyhedral surface is \emph{weakly strictly convex} if for every vertex
$\p_ i$ there is a (support) plane that intersects the surface at exactly
$\p_ i$.  If it is also true that every edge $e$ of the triangulated surface has
a (support) plane that intersects the surface at exactly $e$, we say the
surface is \emph{strongly strictly convex}.  If there is an edge such that
the internal dihedral angle is greater than $180^{\circ}$, we say that edge is
a \emph{non-convex edge} of the surface. 

In addition to being rigid, strongly strictly convex polyhedra with all faces triangles are {\it
  infinitesimally rigid}: 
there is no non-trivial first-order deformation that is  an infinitesimal
  congruence on each triangular face. This point, which was first proved by
  Dehn 
\cite{dehn-konvexer}, is important  in Alexandrov's
subsequent theory concerning the induced metrics on convex polyhedra
(and from there on convex bodies, see \cite{alex}).  Alexandrov also showed
  that Dehn's Theorem can be extended to the case when the polyhedral surface
  is weakly strictly convex, as well as being convex.  In other words,
  vertices of the subdivision can only be vertices of the convex set, and they
  cannot appear in the interior of faces, for example.  If vertices of a
  convex polyhedral surface do lie in the interior of a face, then the surface
  is rigid, but not infinitesimally rigid.  This shows that the underlying
  framework is what determines infinitesimal rigidity, rather than simply the
  surface as a space.   

Our main motivation here is a question concerning the infinitesimal
rigidity of a class of frameworks determined by polyhedra which are weakly strictly convex.

\begin{qu} \label{q:main}
Let $P\subset \E^3$ be a polyhedral surface, with vertices $\p_ 1, \cdots, \p_
n$, such 
that: 
\begin{enumerate}[i.)]
\item \label{weak} $P$ is weakly convex;
\item $P$ is \emph{decomposable}, i.e., it can be written as the union of 
non-overlapping convex polyhedra, such that any two intersect in a common
face, without adding any new vertices. 
\end{enumerate}
Is $P$ then necessarily infinitesimally rigid ? 
\end{qu}

This question comes from \cite{rcnp}, where it is proved that the answer
is positive if the condition, that
there exists an ellipsoid which contains no vertex of $P$ but intersects
all its edges, is added. The goal pursued here is to prove that the answer is
also positive for two classes of polyhedra which are by construction
decomposable. 

\paragraph{Denting polyhedra.}

There is an easy way to construct many examples of polyhedra for which
Condition \ref{weak}.) above holds: start from a convex polyhedron in $\E^3$ and
``dent'' it at some of its edges, in the following manner.

\begin{defi}
Let $P\subset \E^3$ be a polyhedron. $P$ is obtained by \emph{denting} a 
convex polyhedron $Q$ at an edge $e$ if $P$ has the same vertices as $Q$,
and the same faces, except that the two faces of $Q$ adjacent to $e$ (which
are required to be triangles) are replaced by the two other triangles,
sharing an edge, so that the union of the two new triangles has the same
boundary as the union of the two triangles which were removed.
\end{defi}

\begin{figure}[ht]
\begin{center}
\centerline{\psfig{figure=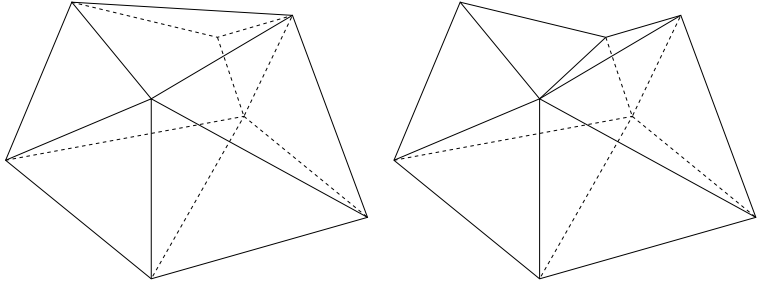,width=10cm}}
\end{center}
\caption{Denting a polyhedron.}
\end{figure}

\paragraph{Simple dented polyhedra are decomposable.}

Clearly, polyhedra obtained by denting a strongly strictly convex polyhedron at any set of edges
have a convex set of vertices (Condition \ref{weak}.) in the question above). Moreover,
those obtained by denting at one edge are decomposable, and this remains true
when denting has occured at two edges which are adjacent to a vertex.

\begin{remark}
Let $P\subset \E^3$ be a convex polyhedron, let $\p_i$ be a vertex of $P$, 
and let $e,e'$ be two edges of $P$ containing $\p_i$ as one of their
endpoints, but which are not two edges of a face of $P$. Let $Q$ be
the polyhedron obtained by denting $P$ at $e$ and at $e'$. Then 
$Q$ is decomposable.
\end{remark}

\begin{proof}
Clearly $Q$ is star-like with respect to $\p_i$, so it can decomposed
as a union of pyramids, each one corresponding to one of the 
faces of $P$ which are not adjacent to $\p_i$.
\end{proof}

However there is no reason to believe that denting a convex polyhedron at
more than two edges, or at two edges not containing a vertex, yields a
decomposable polyhedron.

\paragraph{Infinitesimal rigidity of simple dented polyhedra.}

The first result of this paper is that the answer to Question \ref{q:main}
is positive for polyhedra constructed in this simple manner.

\begin{thm} \label{tm:main}
Let $Q\subset \E^3$ be a polyhedron obtained by denting a strongly convex triangulated polyhedral surface at
one edge, or at two edges sharing a vertex (but which are not both 
contained in a face). Then $Q$ is infinitesimally rigid. 
\end{thm}

The proof, given in Section 2, uses all the excess strength of the
Legendre-Cauchy argument. 

\paragraph{Suspensions.}

A \emph{suspension} $P$ is a polyhedral surface obtained from a closed
polygonal curve $(\p_ 1, \dots, \p_ n)$ connected to two vertices $\mb{N}$, the north
pole, and $\mb{S}$, the south pole.  The edges of $P$ are $[\p_ i, \p_ {i+1}]$, where
$i$ is taken $\bmod$ $n$, and $[\mb{N}, \p_ i], [\mb{S},\p_ i]$ for $i=1, \dots, n$.

Our second task here is to verify Question \ref{q:main} 
in the case of a suspension,
where the natural choice of a decomposition is taken, namely by using the
decomposition  $[\mb{N}, \mb{S}, \p_ i,\p_ {i+1}] $, for $i=1, \dots, n$. 

\begin{thm} \label{tm:suspension}
  Let $P$ be a weakly strictly convex suspension such that the segment $[\mb{N},\mb{S}]$ is contained in $P$. Then $P$ is infinitesimally
  rigid.
\end{thm}

The proof, which is given in Section 3, uses the notion of tensegrity.

\begin{cor}
  Let $P$ be a polyhedron with a convex set of vertices. Suppose that
  $P$ can be cut into simplices (with disjoint interior) so that only
  one new edge, interior to $P$, appears. Then $P$ is infinitesimally
rigid. 
\end{cor}

\begin{proof}
Let $e$ be the edge interior to $P$ which appears in the simplicial
decomposition. Let $P'$ be the union of the simplices, appearing in the
simplicial decomposition of $P$, which contain $e$.
Any isometric first-order deformation of $P$ clearly restricts to
an isometric first-order deformation of $P'$. But $P'$ is by 
construction a suspension, to which Theorem \ref{tm:suspension}
applies. So the isometric first-order deformation of $P'$ is 
trivial, and so is the isometric first-order deformation fo $P$.
\end{proof}

\paragraph{Another type of argument.}

In Section 4 we give another proof of a special case of Theorem
\ref{tm:suspension}, based on a simple idea: given a suspension 
we compute the first order variation, under a variation of the $\mb{N}-\mb{S}$
distance, of the sum of the angles of the simplices at the $\mb{N}-\mb{S}$ line.
This leads us in Section 5 to define a symmetric matrix attached 
to a polyhedron along with a simplicial decomposition, which is 
non-singular if and only if the polyhedron is infinitesimally 
rigid. We then formulate another question, for which a positive answer
would imply a positive answer to Question \ref{q:main}. 

\section{Rigidity of dented polyhedra}

\paragraph{Outline.}

The proof of Theorem \ref{tm:main} follows quite precisely the arguments
of the original Cauchy proof, but with slightly sharper estimates. 
We consider a first-order deformation
of $Q$, and associated to each edge $e$ of $Q$ a ``sign'', which is 
$0$ if the dihedral angle at $e$ does not vary (at first order), 
$+$ if it increases, $-$ if it decreases. The first lemma is of 
a geometrical nature: 

\begin{lemma} \label{lm:geom} The following hold:
\begin{enumerate}
\item For each vertex $\q_i$ of $Q$ where $Q$ is not convex, either all the signs
  attached to the edges of $Q$ containing $\q_i$ are $0$, or some have a
  sign $+$ and others have a sign $-$. In that case there are at least 
  3 edges containing $\q_i$ with a sign which is not $0$.
\item At each vertex $\q_i$ of $Q$ where $Q$ is convex, either the signs
  assigned to 
all edges are $0$, or there are at least 4 changes of signs when one considers
the edges containing $\q_i$ in cyclic order.
\end{enumerate}
\end{lemma}

Here a ``sign change'' means a sequence of edges such that the sign on the
first edge is $+$ (resp. $-$), the sign on the last edge is $-$ (resp. $+$),
and the signs on the other edges are all $0$.

The second lemma is of a topological nature. 

\begin{lemma} \label{lm:topo}
Let $\Gamma$ be a graph embedded in the sphere, which is the 1-skeleton
of a cellular decomposition of $\mb{S}^2$. It is not possible to assign a sign 
$+$ or $-$ to each edge of $\Gamma$ such that:
\begin{itemize}
\item the signs assigned to all the edges containing a given vertex are never
  all the same,
\item there are at least 4 changes of sign at all vertices except at most
  three.
\end{itemize}
\end{lemma}

\paragraph{Proof of the geometric lemma.}

The geometric lemma \ref{lm:geom} follows from a remark concerning
infinitesimal deformations of spherical polygon (see e.g. 
\cite{gluck-generique,polygones}).

\begin{lemma} \label{lm:poly}
Let $\p\subset \mb{S}^2$ be a polygon, with vertices $\p_ 1, \cdots, \p_ n$. Let
$\theta_1, \cdots, \theta_n$ be the angle of $\p$ at the vertices. The
possible first-order variations of those angles under the infinitesimal
deformations of $\p$ are the $n$-uples $(\theta'_1, \cdots, \theta'_n)$
characterized by the relation: 
$$ \sum_{i=1}^n \theta'_i \p_ i =0~. $$
\end{lemma}

Lemma \ref{lm:geom} now follows by considering the link of each
vertex of $Q$. If $\q_i$ is such a vertex, then, since the
set of vertices of $Q$ is convex, the link $L(\q_i)$ of $Q$ at $\q_i$ is 
contained in an open hemisphere, so it is clear from Lemma 
\ref{lm:poly} that the signs associated to the edges of $Q$
containing $\q_i$ --- which correspond to the vertices of $L(\q_i)$,
with a sign $+$ if the angle of $L(\q_i)$ increases and $-$ if
it decreases --- cannot be all non-positive or all non-negative. Moreover, if
$Q$ is convex at $\q_i$, then there are at least 4 changes of
signs; otherwise there would be exactly 2 changes of signs,
and, if $\mb{u} \in \R^3$ were a vector such that the plane orthogonal
to $\mb{u}$ containing $\q_i$ separated the edges with a $+$ from the
edges with a $-$, then the scalar product of $\mb{u}$ with the
left-hand side of the equation in Lemma \ref{lm:poly} would
be non-zero.

It is interesting to observe that Lemma \ref{lm:topo} also easily follows from an argument using stresses as in Section \ref{Sect:suspensions}.  

\paragraph{Proof of the combinatorial lemma.}

The proof of Lemma \ref{lm:topo} follows quite directly the
argument given by Cauchy. Let $v,e,f$ be the number of vertices,
edges and faces of $\Gamma$, respectively, and let $s$ be the number
of changes of signs, i.e., pairs $\{ e,e'\}$ of edges which are adjacent
edges of a face, with a sign $+$ on one and $-$ on the other. 

By the hypothesis of the Lemma, there are at least 4 changes of
signs at each vertex except perhaps at 3 vertices where there are
at least 2 changes of signs, so that:
$$ 4v-6 \leq s~. $$

However there is an upper bound on the number of changes of
signs on the faces of $\Gamma$: there can be at most 2 changes
of signs on a triangular face, at most 4 changes of signs on 
a face with 4 or 5 edges, at most 6 changes of signs on a face
with 6 or 7 edges, etc. Calling $f_k$ the number of faces with
$k$ edges, this means that:
$$ s \leq 2f_3 + 4f_4 + 4 f_5 + 6f_6 + 6f_7 + 8f_8 + 8f_9+ \cdots~. $$
However each edge of $\Gamma$ bounds two faces, which shows that:
$$ 2e = 3f_3 + 4f_4 + 5f_5 + 6f_6 + 7f_7 + \cdots~. $$
Taking twice this equation and substracting:
$$ 4f = 4f_3 + 4 f_3 + 4f_4 + 4f_5 + 4f_6 + 4f_7 + \cdots~, $$
we obtain that:
$$ 4e-4f = 2f_3 + 4f_4 + 6f_5 + 8f_6 + 10 f_7 + \cdots~, $$
so that:
$$ 4e-4f\geq  s~. $$
Putting together the two inequalities on $s$ yields that 
$4e-4f\geq 4v-6$, which contradicts the Euler relation, 
$v-e+f=2$.

\paragraph{Proof of the rigidity theorem.}

The proof of Theorem \ref{tm:main} now follows as in the original
Cauchy proof. Suppose that $Q$ has a non-trivial infinitesimal 
deformation, and assign a sign $+$, $-$ or $0$ to each edge
depending on whether the angle at that edge increases, decreases
or stays constant in the deformation. Then consider the graph $\Gamma$
obtained from the 1-skeleton of $Q$ by removing all edges with a $0$
and all vertices of $Q$ which are contained only in edges with sign $0$.
It follows from Lemma \ref{lm:poly} that $\Gamma$ is still the
1-skeleton of a cell decomposition of the sphere, because any vertex
of $\Gamma$ is contained in at least 3 edges of $\Gamma$.
Lemma \ref{lm:geom} shows that there are at least 4 changes of sign
at each vertex of $\Gamma$ except perhaps at 3 vertices where there
are at least 2 changes of sign, and Lemma \ref{lm:topo} shows that
this is impossible.

\section{Suspensions}\label{Sect:suspensions}

We introduce the definition of infinitesimal rigidity in the context of a tensegrity.  Consider a finite collection points $\p=(\p_  1, \dots, \p_  n)$ in $\E^d$ and a graph $G$ with those points as vertices, and with edges, called bars, cables or struts,  between some pairs of those points. Consider vectors $\p'=(\p'_ 1, \dots, \p'_ n)$ in $\E^d$,  $\p'_ i$  regarded as an infinitesimal motion, a velocity, associated to $\p_  i$ for each $i = 1, \dots, n$.  We say that $\p'$ is an \emph{infinitesimal flex} of the tensegrity if the following equation of vector inner products holds:
\begin{equation}\label{infdef}
(\p_  i - \p_  j)\cdot (\p'_ i - \p'_ j)  \begin{cases} 
 \le 0 & \text{if $\{i,j\}$ is a cable,}\\
 = 0 & \text{if $\{i,j\}$ is a bar,}\\
 \ge 0 & \text{if $\{i,j\}$ is a strut.}\\
 \end{cases}
\end{equation}
A bar framework is \emph{infinitesimally rigid} if the only infinitesimal flexes are the trivial ones that come as the derivative of a family of rigid congruences of all of Euclidean space restricted to the configuration $\p$.  

One useful tool in showing infinitesimal rigidity is the concept of a
\emph{stress}, which are scalars $\omega_{ij}$ associated to each edge
$\{i,j\}$.  We write this as a single vector  $\omega = (\dots, \omega_{ij},
\dots)$.    We say $\omega$ is an \emph{equilibrium stress} if the following
vector equation holds for each $i$: 
\begin{equation}\label{equilibrium}
\sum_j \omega_{ij} (\p_ i - \p_ j) = 0 \quad  \text{for every $\{i,j\}$ an edge of
  $G$.} 
\end{equation}
We say the the stress $\omega = (\dots, \omega_{ij}, \dots)$ is \emph{proper}
if $\omega_{ij}$ is non-negative for cables and non-positive for struts.
(There is no sign condition for bars.) 

One important way to use this concept is the following:  

\begin{lemma}\label{block} If $\p'$ is infinitesimal flex of a tensegrity
  framework $G(\p)$ and $\omega = (\dots, \omega_{ij}, \dots)$ is a proper
  equilibrium stress for $\p$, then $\sum_{ij} \omega_{ij} (\p_ i - \p_ j)\cdot
  (\p'_ i - \p'_ j) = 0 $ and thus if  $\omega_{ij} \ne 0$, then $(\p_ i - \p_ j)\cdot
  (\p'_ i - \p'_ j) = 0$. 
\end{lemma}

\begin{proof} The equation  $\sum_{ij} \omega_{ij} (\p_ i - \p_ j)\cdot (\p'_ i -
  \p'_ j) = 0 $ follows by taking the inner product of Equation
  (\ref{equilibrium}) with $\p'_ i$ and summing over all $i$.  The last part of
  the conclusion follows from the condition that $\omega$ is proper and the
  inequality in Condition (\ref{infdef}). 
\end{proof}

One way to look at Lemma \ref{block} is that a proper equilibrium stress in a
tensegrity ``blocks" the infinitesimal motion that would decrease a cable or
increase a strut in the sense that the inequalities of Condition
$(\ref{infdef}$) would be strict.

\begin{cor} \label{exchange} 
  Suppose that a tensegrity has a proper
  equilibrium stress and is infinitesimally rigid.  Then the bar framework
  obtained by removing any of the edges with a non-zero stress, and converting
  all the other edges to bars,  is infinitesimally rigid.
\end{cor}

\begin{proof} 
Any infinitesimal flex that is non-zero on the removed edge is blocked by the
stress.  So the removed edge can be put back as far as the infinitesimal flex
is concerned. 
\end{proof}

One consequence of Corollary \ref{exchange} is that  one bar framework can be
exchanged for another, by adding a bar and removing another when there is an
equilibrium stress that is non-zero on both bars. 

Suppose that we have a suspension $G(\p)$ regarded as polyhedral surface.  We say that a suspension is
\emph{$\mb{N}-\mb{S}$ decomposable} if the projection on the plane 
orthogonal to the line through $\mb{N}$
and $\mb{S}$ of the equator is one-to-one and the projection of the point $\mb{N}$ (and
$\mb{S}$) lies inside the projection of the equator.  So an $\mb{N}-\mb{S}$ decomposable
suspension can be decomposed into non-overlaping tetrahedra, as in the Main
Question \ref{q:main}. 

For any suspension we create a tensegrity by labeling the $\{\mb{N}, \mb{S}\}$ edge and
the equatorial edges as  cables.  The lateral edges are simply bars.  Call
this a \emph{tensegrity suspension}. 

\begin{lemma} \label{convex case} A strongly strictly convex tensegrity
  suspension has a proper equilibrium stress. 
\end{lemma}

\begin{proof} 
  A convex suspension with $n$ vertices (including $\mb{N}$ and $\mb{S}$)
  has $3(n-2)=3n-6$ edges, and the associated tensegrity has $3n-6+1=3n-5$
  edges.  The equilibrium conditions involve $3n-5$ variables and $3n$ linear
  equations, one for each coordinate of each vertex. There is always a
  $6$-dimensional linear subspace of the $3n$-dimensional space that is
  orthogonal to the space of all possible stresses.  (This corresponds to the
  trivial infinitesimal flexes.)  So there must be a one-dimensional space of
  equilibrium stresses.  The signs of the stresses are proper by the
  Cauchy-Dehn argument. There are only four edges incident to an equatorial
  vertex, and so they alternate in sign.  Note also, in this case, that the
  sign of the stresses on the lateral edges is opposite from the sign on the
  stresses of the equator and $\mb{N}-\mb{S}$ edge.  
\end{proof}

\begin{lemma} \label{supstress} Any $\mb{N}-\mb{S}$ decomposable strictly weakly convex
  tensegrity suspension $G(\p)$ has a proper equilibrium stress. 
\end{lemma}

\begin{figure}[ht]
\centerline{\psfig{figure=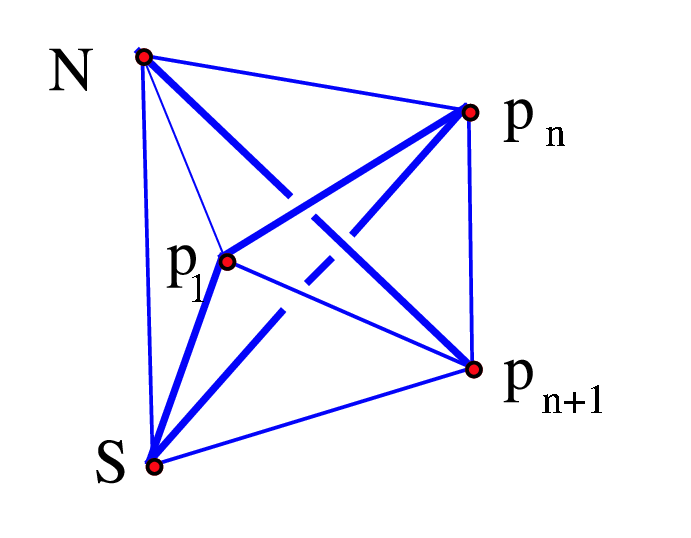,width=6cm}}
\caption{This shows the case when the edge $[\mb{N}, \p_ {n+1}]$ is non-convex.
   The thick edges are struts with a negative stress and the thin edges are
   cables with a positive stress.} 
 \label{fig:viennafig}
\end{figure}

\begin{proof} 
The proof is by induction on $n$, the number of equatorial
vertices.  If there are only $3$ vertices on the equator, the decomposability
condition implies that the suspension is strongly strictly convex, so Lemma
\ref{convex case} implies that it has a proper equilibrium stresss.
Whenever there is a non-convex lateral edge with $n+1$ equatorial edges,  we
will show how to create the proper equilibrium stress from another  $\mb{N}-\mb{S}$
decomposable strictly weakly convex tensegrity suspension with $n$
equatorial edges. 

We assume that the tensegrity suspension on the vertices $\p_ 1, \dots, \p_ n$ has
a proper equilibrium stress, and we wish to show that the tensegrity
suspension on the vertices $\p_ 1, \dots, \p_ n, \p_ {n+1}$ also has a proper
equilibrium stress, where, say, the lateral edge $[\mb{N}, \p_ {n+1}]$ is not
strictly convex. 

Consider the wedge $W$ determined by the two planes $[\mb{N}, \mb{S}, \p_ 1]$ and $[\mb{N}, \mb{S},
\p_ n]$ between $\p_ 1$ and $\p_ n$ cyclicly around the $\mb{N}-\mb{S}$ axis.  Since the
suspension is weakly strictly convex at $\p_ {n+1}$, this point 
can not be in the
tetrahedron determined by $[\mb{N}, \mb{S}, \p_ 1, \p_ n]$, yet it has to be 
in $W$.  So the five
points $\mb{N}, \mb{S}, \p_ 1,\p_ {n+1}, \p_ n$ determine a (small) tensegrity suspension over
a triangle, where  $[\mb{N}, \p_ {n+1}]$ is the axis since it is non-convex.  In this
case $[\mb{N},\mb{S}]$ is a lateral edge, and $[ \p_ 1,\p_ {n+1}]$ is an 
equatorial edge.  Thus
we can choose an equilibrium stress for the small suspension such that the
stress on $[\p_ 1,\p_ {n+1}]$ is exactly the negative of the equilibrium stress
of the suspension on $\p_ 1, \dots, \p_ n$ on $[ \p_ 1,\p_ {n+1}]$.  When these two
stresses are added, they cancel.  So the stresses on $[\mb{N},\mb{S}]$ for both the
small suspension and the large one are positive, and hence their sum is
positive.  So this sum is a proper equilibrium stress for the larger
suspension as desired.  Figure \ref{fig:viennafig} shows this situation.   

If the four vertices $\mb{N}, \p_ 1, \p_ {n+1}, \p_ n$ are coplanar, then there is an
equilibrium stress as before, except that it is $0$ on $[\mb{N}, \mb{S}]$  and the same
argument applies. 

If there are no non-convex lateral edges in the suspension, since the
equatorial edges are always convex,  the whole suspension is convex and Lemma
\ref{convex case} applies.   
\end{proof}

\begin{cor}  Any $\mb{N}-\mb{S}$ decomposable strictly weakly convex suspension is
   infinitesimally rigid. 
\end{cor}

\begin{proof} 
  Since the bar framework obtained by adding the $[\mb{N},\mb{S}]$ bar is
  infinitesimally rigid and there is an equilibrium stress non-zero on
  $[\mb{N},\mb{S}]$ when that bar is added, Corollary \ref{exchange} implies that
  the suspension itself is infinitesimally rigid. 
\end{proof}

\section{Rigidity of suspensions through variations of angles}

\paragraph{Suspensions over convex polygons.}

Consider a suspension $P$ with vertices $\mb{N},\mb{S}$ and $\p_ 1, \cdots, \p_ n$, 
as defined above. There is a natural subclass of suspensions, which 
are all weakly convex and decomposable.

\begin{defi}
$P$ is a {\bf suspension over a convex polygon} if the projection of the
closed polygonal line with  
vertices $\p_ 1, \cdots, \p_ n$ (in this cyclic order), along 
$(\mb{N},\mb{S})$, on any plane transverse
to the line $(\mb{N},\mb{S})$ is convex and contains in its interior the
intersection of the plane with $(\mb{N},\mb{S})$.
\end{defi}

\begin{thm} \label{tm:elementary}
Let $P$ be a suspension over a convex polygon. Then $P$ is infinitesimally rigid.
\end{thm}

Note that, under the hypothesis of the theorem, $P$ is always decomposable
(because it is a suspension) and it is also weakly convex. However,  the
hypothesis in Theorem \ref{tm:elementary} is stronger than in Theorem
\ref{tm:suspension} since the projection of the equator on 
any plane transverse to the $\mb{N}-\mb{S}$ axis 
is required to be a convex polygon.
So the statement of this theorem is not very interesting in itself,
we include it here because its proof is different from the proof given
for Theorem \ref{tm:elementary} an can be interesting as an indication
of possible ways to tackle Question \ref{q:main} (as seen in the next
section). 

The first step in the proof is to apply a projective transformation
to $P$ so that the planes orthogonal to $[\mb{N},\mb{S}]$ at $\mb{N}$ and at $\mb{S}$ are
support planes of $P$. This does not change the infinitesimal rigidity
or flexibility of $P$ since it is well known (at least since works of
Darboux \cite{darboux12},  Sauer \cite{sauer}, and J. Clerk Maxwell in the 19-th Century, that infinitesimal rigidity
is a projectively invariant property). From here on we suppose that
this additional property is satisfied.

\paragraph{Deformations of simplices.}

The argument given below is based on a computation concerning
the first-order deformations of Euclidean simplices. We consider 
a simplex with two vertices called $\mb{N}$ and $\mb{S}$ of coordinates
$(0,0,1)$ and $(0,0,0)$ in $\R^3$, and with two other vertices $\p_1$
and $\p_2$ of coordinates $(r_1\cos(\alpha_1),r_1\sin(\alpha_1),z_1)$
and $(r_2\cos(\alpha_2),r_2\sin(\alpha_2),z_2)$.

\begin{lemma} \label{lm:simplex}
There is a unique first-order deformation of this simplex under which the
distance between $\mb{N}$ and $\mb{S}$ varies at speed $1$, and the lengths of all
the other edges remain constant. Under this deformation, the first-order
variation of the angle $\theta=\alpha_2-\alpha_1$ at the edge $\mb{N}-\mb{S}$ is:
$$ \theta' = \frac{1}{r_1r_2\sin\theta}\left(
(z_1-z_2)^2 + z_1(1-z_1)\left(1-\frac{r_2}{r_1}\cos\theta\right) +
z_2(1-z_2)\left(1-\frac{r_1}{r_2}\cos\theta\right)
\right)~. $$
\end{lemma}

\begin{proof}
The square of the distance between $\mb{S}$ and $\p_1$ is equal to $z_1^2 +
r_1^2$. Since this remains constant under the deformation we have: 
$$ z_1 z'_1 + r_1 r'_1 = 0~. $$
Similarly, calling $z$ the third coordinate of $\mb{N}$ (so that $z=1$ 
``before'' the deformation takes place), 
the square of the distance between $\mb{N}$ and $\p_1$ is equal 
to $(z-z_1)^2 + r_1^2$, and it is constant under the deformation, so that:
$$ r_1 r'_1 + (1-z_1)(1-z'_1) = 0~. $$
It follows that
$$ z'_1 = 1-z_1~, ~~ r'_1 = - \frac{z_1}{r_1} (1-z_1)~, $$
and similarly:
$$ z'_2 = 1-z_2~, ~~ r'_2 = - \frac{z_2}{r_2} (1-z_2)~. $$

Furthermore the square of the distance between $\p_1$ and $\p_2$ is equal to
$(z_2-z_1)^2+r_1^2+r_2^2-2r_1r_2\cos\theta$. Since it remains constant in the
deformation we have that:
$$ (z_2-z_1)(z'_2-z'_1) + r_1r'_1 + r_2r'_2 -r'_1r_2\cos\theta -
r_1r'_2\cos\theta + r_1r_2\theta' \sin\theta=0~, $$
so that:
$$ -(z_2-z_2)^2 +\left(1-\frac{r_2}{r_1}\cos\theta\right) r_1r'_1 
+ \left(1-\frac{r_1}{r_2}\cos\theta\right) r_2r'_2 +  r_1r_2\theta'
\sin\theta=0~, $$
and thus:
$$ r_1r_2\theta' \sin\theta = (z_2-z_2)^2 +
z_1(1-z_1)\left(1-\frac{r_2}{r_1}\cos\theta\right) + 
z_2(1-z_2)\left(1-\frac{r_1}{r_2}\cos\theta\right)~, $$
from which the result follows.
\end{proof}

\paragraph{An invariant controling the infinitesimal rigidity  of
  suspensions.} 

Let $\p_1, \cdots, \p_n$ be the vertices of $P$ different from $\mb{N}$ 
and $\mb{S}$, in the cyclic order on which they appear on the ``equator'' 
of $P$. Suppose that the coordinates of $\p_ i$ are $(r_i\cos\alpha_i, r_i\sin
\alpha_i, z_i)$ for $1\leq i\leq n$. Below we use cyclic notation, 
so that $\p_{n+1}=\p_1$.

We define a quantity $\Lambda(P)$ which is the sum of the terms appearing
in Lemma \ref{lm:simplex} above for the simplices $[\mb{N},\mb{S},\p_i, \p_{i+1}]$. 
We will see that $\Lambda(P)=0$
if and only if $P$ is infinitesimally flexible. Below we will also give
simpler geometric expressions of $\Lambda$, leading in particular to the
proof of Theorem \ref{tm:elementary}. It is defined as:
$$ \Lambda(P) := \sum_{i=1}^n \frac{1}{r_i r_{i+1} \sin\theta_i} 
\left(
(z_{i+1}-z_i)^2 + z_i(1-z_i)\left(1-\frac{r_{i+1}}{r_i}\cos\theta_i\right) +
z_{i+1}(1-z_{i+1})\left(1-\frac{r_i}{r_{i+1}}\cos\theta_i\right)\right)~. $$

Consider any first-order deformation of $P$. If the distance between
$\mb{N}$ and $\mb{S}$ does not vary, then the deformation is trivial, because
each simplex $S_i$ of vertices $\mb{N},\mb{S},\p_i$ and $\p_{i+1}$  would then remain 
the same. It follows that the angle $\theta_i$ of $S_i$ at the edge
$\mb{N}-\mb{S}$ varies accordingly to Lemma \ref{lm:simplex}. Since the sum
of the angles of the simplices $S_i$ at the edge $\mb{N}-\mb{S}$ has to 
be equal to $2\pi$, a first-order deformation of $P$ is trivial unless
the sum of the terms corresponding to Lemma  \ref{lm:simplex} for 
the simplices $S_i, 1\leq i\leq n$, sum to $0$. It is not difficult
to check that the converse is true, too: if the sum of the first-order
variations of the angle at $\mb{N}-\mb{S}$ of the simplices $S_i$ vanishes, 
then there is non-trivial first-order deformation of $P$. This shows
the following statement:

\begin{lemma} \label{lm:lambda}
$P$ is infinitesimally rigid if and only if $\Lambda(P)\neq 0$.
\end{lemma}

\paragraph{Another expression of $\Lambda(P)$.}

Now consider the polygon $\p$ which is the orthogonal 
projection of the ``equator'' of $P$ on the plane $z=0$. 
Let $\mb{u}_1, \cdots, \mb{u}_n$ be its
vertices, with $ \mb{u}_i$ equal to the orthogonal projection of $\p_i$ on $z=0$.
We call $a_i$ twice the area of the triangle $( \mb{0}, \mb{u}_i,\mb{u}_{i+1})$, and $b_i$
twice the {\it oriented} area of the triangle $(\mb{u}_{i-1},  \mb{u}_i, \mb{u}_{i+1})$. 

\begin{lemma} \label{lm:expression}
We have:
$$ \Lambda(P) = \sum_{i=1}^n \frac{(z_{i+1}-z_i)^2}{a_i} +
\frac{z_i(1-z_i)b_i}{a_{i-1} a_i}~. $$
\end{lemma}

\begin{proof}
A simple computation shows that:
\begin{eqnarray*}
\Lambda(P)
& = & \sum_{i=1}^n \frac{(z_{i+1}-z_i)^2}{r_i r_{i+1} \sin\theta_i} +
\frac{z_i(1-z_i)}{r_i^2} \left(
\frac{r_i-r_{i+1}\cos\theta_i}{r_{i+1}\sin\theta_i} + 
\frac{r_i-r_{i-1}\cos\theta_{i-1}}{r_{i-1}\sin\theta_{i-1}}
\right) \\
& = & \sum_{i=1}^n \frac{(z_{i+1}-z_i)^2}{a_i} +
\frac{z_i(1-z_i)}{r_i^2} \left(
\frac{r_i-r_{i+1}\cos\theta_i}{r_{i+1}\sin\theta_i} + 
\frac{r_i-r_{i-1}\cos\theta_{i-1}}{r_{i-1}\sin\theta_{i-1}}
\right)~. 
\end{eqnarray*}

Let $\beta_i$ be the angle at $\mb{u}_i$ between $(\mb{u}_i,\mb{u}_{i-1})$ and $(\mb{u}_i,\mb{0})$
and let $\gamma_i$ be the angle at $\mb{u}_i$ between the oriented lines $(\mb{u}_i,0)$
and $(\mb{u}_i,\mb{u}_{i+1})$.
It is easy to check on a diagram
that $r_{i+1}\sin\theta_i/(r_i-r_{i+1}\cos\theta_i)$ is the tangent
of $\gamma_i$, while $r_{i-1}\sin\theta_{i-1}/(r_i-r_{i-1}\cos\theta_{i-1})$
is the tangent of $\beta_i$. It follows that:
\begin{eqnarray*}
\Lambda(P) & = & 
\sum_{i=1}^n \frac{(z_{i+1}-z_i)^2}{a_i} +
\frac{z_i(1-z_i)}{r_i^2} (\cot(\gamma_i)+\cot(\beta_i)) \\
& = & \sum_{i=1}^n \frac{(z_{i+1}-z_i)^2}{a_i} +
\frac{z_i(1-z_i) \sin(\alpha_i)}{r_i^2 \sin(\beta_i)\sin(\gamma_i)}~, 
\end{eqnarray*}
where $\alpha_i$ is the interior angle of $\p$ at $\mb{u}_i$, i.e.,
$\alpha_i = \beta_i + \gamma_i$. 
Note that $a_i=r_i \|\mb{u}_{i+1}-\mb{u}_i\| \sin\theta_i$ and that 
$a_{i-1} = r_{i-1} \| u_i-u_{i-1}\| \sin\theta_{i-1}$, it follows that
$$ \Lambda(P) = \sum_{i=1}^n \frac{(z_{i+1}-z_i)^2}{a_i} +
\frac{z_i(1-z_i) \sin(\alpha_i)\|\mb{u}_{i+1}-\mb{u}_i\|.
\| \mb{u}_i-\mb{u}_{i-1}\|}{a_{i-1} a_i}~, $$
and the expression given in the Lemma follows because $b_i=\sin(\alpha_i)\|
\mb{u}_{i+1}-\mb{u}_i\|.\| \mb{u}_i-\mb{u}_{i-1}\|$. 
\end{proof}

The proof of Theorem \ref{tm:elementary} is a direct consequence of 
Lemma \ref{lm:expression} since, under the hypothesis of the theorem,
all terms involved in the sum are non-negative (and some of them
are strictly positive).

\section{A more technical version of Question \ref{q:main}}

The arguments in the previous section suggest a way to transform 
Question \ref{q:main} in a way which makes it much less elementary
but perhaps more precise in a fairly technical manner. For a 
general decomposable polyhedron $P$, along with a decomposition as a
union of simplices with disjoint interior (and with no vertex 
beyond those of $P$) it is still possible to define an invariant
generalizing the number $\Lambda$ defined above for suspensions. 
However in this more general case $\Lambda(P)$ is a $r\times r$
matrix, where $r$ is the number of ``interior'' edges in the 
simplicial decomposition of $P$. 

After defining $\Lambda(P)$,
we will show here that it is a symmetric matrix, and that $P$ is
infinitesimally rigid if and only if 
$\Lambda(P)$ is non-singular. In addition, the content of the previous
section shows that, under some mild assumptions,
its diagonal has positive entries. This suggests that $\Lambda(P)$ might
be positive definite whenever $P$ is weakly convex; this would at least
imply a positive answer to Question \ref{q:main}.

\paragraph{An extension of the invariant $\Lambda$.}

In this section we consider a polyhedron $P$, along with a simplicial
decomposition of $P$ with no vertex except the vertices of $P$. 
We call $S_1, \cdots, S_q$ the simplices in the decomposition of 
$P$, and $e_1, \cdots, e_r$ the interior edges
of the decomposition, i.e., the segments which are edges of the 
simplicial decomposition but are contained in the interior of $P$.

Clearly any isometric first-order deformation of $P$ is uniquely determined by
the first-order variation of the lengths of the $e_i$, because the other edges
of the simplicial decomposition have fixed length. 

\begin{defi}
  Let $E\subset \R^r$ be the set of lengths $l_1,\cdots, l_r$ so that,
for each of the simplices $S_i$, if the lengths of the edges of $S_i$ which
are in the interior of $P$ are fixed at the values determined by the $l_j$ 
and the length of the other edges are equal to their length in $P$, the
resulting 6 numbers are indeed the lengths of the edges of an Euclidean
simplex. 
\end{defi}

There is a ``special'' element $l^0=(l_1^0,\cdots, l_r^0)$ in $E$,
it is the $r$-uple of lengths of the edges $e_i$ in the polyhedron
$P$.

To each choice of $(l_1, \cdots, l_r)\in E$ we can associate an Euclidean 
structure on the simplicial decomposition of $P$ used here, but with
cone singularities at the edges $e_i$. So to each interior edge $e_i$
is attached a number, the sum of the angles at $e_i$ of the simplices
containing it (or in other terms the angle around the cone singularity
corresponding to $e_i$), which we call $\theta_i$. This defines a 
map:
$$ 
\begin{array}{rrcl}
\phi:& E &\mapsto & (0,\infty)^r \\
& (l_1, \cdots, l_r)& \rightarrow & (\theta_1, \cdots, \theta_r)
\end{array} $$

It clearly follows from the construction that $\phi(l^0)=(2\pi, \cdots,
2\pi)$. 

\begin{defi}
Let $\Lambda(P)$ be the Jacobian matrix of $\phi$ at $l^0$, i.e.:
$$ \Lambda(P) := \left( \frac{\partial \theta_i}{\partial l_j}\right)_{1\leq
  i,j\leq r}~. $$
\end{defi}

In the special case of suspensions considered in Section 3 and 4, $\Lambda(P)$
is a $1\times 1$ matrix, and its unique entry is the quantity called
$\Lambda(P)$ there. 

\paragraph{$\Lambda(P)$ and the infinitesimal rigidity of $P$.}

As for suspensions, $\Lambda(P)$ can be used to determine when $P$ is
infinitesimally rigid. 

\begin{lemma}
$P$ is infinitesimally flexible if and only if $\Lambda(P)$ is singular (i.e.,
its kernel has dimension at least $1$).  
\end{lemma}

\begin{proof}
Suppose first that $P$ is not infinitesimally rigid, and consider a 
 first-order edge-length preserving deformation. Let $(l'_1, \cdots, l'_r)$ be the
corresponding first-order variations of the lengths of the interior
edges of the simplicial decomposition. Under the same deformation,
the sum of the angles at each interior edge $e_i$ remains equal
to $2\pi$ (at first order) so that $\Lambda(P)(l'_1, \cdots, l'_r)=0$,
and this shows that $\Lambda(P)$ is singular.

Suppose conversely that $\Lambda(P)$ is singular, and let $(l'_1, \cdots,
l'_r)$ be a non-zero element in the kernel of $\Lambda(P)$. It defines
a first-order variation of the Euclidean metric on each of the simplices
appearing in the decomposition of $P$, and therefore of the Euclidean
structure, with singularities at the edges, naturally defined on this
simplicial complex. However the first-order variation under $(l'_1, \cdots,
l'_r)$ of the angle
at each of the interior edges vanishes precisely because $(l'_1, \cdots,
l'_r)$ is in the kernel of $\Lambda(P)$. This means that the 
Euclidean structure remains associated to an Euclidean polyhedron,
and therefore that $P$ is not infinitesimally rigid.
\end{proof}

\paragraph{$\Lambda(P)$ is symmetric.}

This is another simple property of $\Lambda(P)$. To prove it we need
some additional notation. We call $\eb_1, \cdots, \eb_{\rb}$ the 
edges of the simplicial decomposition which are contained in the
boundary of $P$ -- i.e., those which are not among the $e_i$ -- and
$\lb_1, \cdots, \lb_{\rb}$ their lengths. For each simplex $S_i$ and
each edge $e_j$ (resp. $\eb_j$) which is an edge of $S_i$ we call
$\alpha_{i,j}$ (resp. $\alphab_{i,j}$) the angle of $S_i$ at $e_j$ (resp. at
$\eb_j$). 

We then introduce the sum of the ``mean curvatures'' of the simplices:
$$ H := \sum_{i,j} l_j\alpha_{i,j} + \sum_{i,j} \lb_j\alphab_{i,j}~, $$
where the first sum is over all simplices $S_i$ and edges $e_j$ such that
$e_j$ is an edge of $S_i$, while the second sum is the corresponding quantity
with $\eb_j$ instead of $e_j$. Then, under a first-order deformation:
$$ dH = \sum_{i,j} (l_jd\alpha_{i,j} + \alpha_{i,j}dl_j) + 
\sum_{i,j} (\lb_jd\alphab_{i,j} + \alphab_{i,j}d\lb_j)~. $$
Now we consider $H$ as a function over $E$, which means that we fix the
values of the $\lb_j$, the lengths of the edges of the simplicial
decomposition of $P$ which are on the boundary of $P$. The formula
for the first-order variation of $H$ simplifies and becomes:
$$ dH = \sum_{i,j} \alpha_{i,j}dl_j + \sum_{i,j}
l_jd\alpha_{i,j} + \sum_{i,j} \lb_jd\alphab_{i,j}~. $$
But the celebrated Schl\"afli formula states that
$$ \sum_{i,j} \lb_jd\alphab_{i,j} + \sum_{i,j}
l_jd\alpha_{i,j} =0~, $$
so that
$$ dH = \sum_{i,j} \alpha_{i,j}dl_j = \sum_j \theta_j dl_j~. $$
This means that $\Lambda(P)$, as defined above, is the Hessian 
matrix of $H$, considered as a function of the $l_j$ (which are 
coordinates on $E$). So it is a symmetric matrix.

\paragraph{Another question.}

It follows from Section 4 that, under some fairly simple geometric
hypothesis on $P$, the diagonal of $\Lambda(P)$ is positive. This
leads to the

\begin{qu}
Is $\Lambda(P)$ positive definite whenever $P$ is weakly convex ? 
\end{qu}

A positive answer would imply that weakly convex and decomposable
polyhedra are infinitesimally rigid, i.e., a positive answer to Question
\ref{q:main}.

\bibliographystyle{alpha}

\def\cprime{$'$}

\end{document}